\documentclass[12pt, reqno]{amsart}

\usepackage{version}
\usepackage[T1]{fontenc} 
\usepackage[utf8]{inputenc} 
\usepackage{amsmath}
\usepackage{amssymb}
\usepackage{mathrsfs}
\usepackage{amsthm}
\usepackage{latexsym}
\usepackage{wasysym}
\usepackage{dcpic, pictex}
\usepackage{graphicx}
\usepackage{tikz-cd}

\newtheorem{theorem}{Theorem}[section]
\newtheorem{lemma}[theorem]{Lemma}
\newtheorem{corollary}[theorem]{Corollary}
\newtheorem{proposition}[theorem]{Proposition}

\theoremstyle{definition}
\newtheorem{definition}[theorem]{Definition}

\newtheorem{remark}[theorem]{Remark}

\newtheorem{conjecture}[theorem]{Conjecture}

\author{Matteo Fiacchi}

\def\C{{\mathbb{C}}}

\def\N{{\mathbb{N}}}
\def\B{{\mathbb{B}}}
\def\Bn{{\mathbb{B}^{n}}}
\def\Cn{{\mathbb{C}^{n}}}
\def\S{{\mathcal{S}}}

\newcommand{\id}{{\sf Id}}

\begin{document}
\title[Embedding and approximation]{The embedding conjecture and the approximation conjecture in higher dimension} 

\address{M. Fiacchi: Dipartimento di Matematica\\
Universit\`{a} di Roma \textquotedblleft Tor Vergata\textquotedblright\ \\Via Della Ricerca Scientifica 1, 00133 \\
Roma, Italy} \email{fiacchi@mat.uniroma2.it}
\begin{abstract}
In this paper we show the equivalence among three conjectures (and related open questions), namely, the embedding of univalent maps of the unit ball into Loewner chains, the approximation of univalent maps with entire univalent maps and the immersion of domain biholomorphic to the ball in a Runge way into  Fatou-Bieberbach domains. 
\end{abstract}
\maketitle

\section{Preliminaries}
Let $\Bn:=\{z\in\Cn: ||z||<1\}$ be the unit ball of $\Cn$ and 
$$\S:=\{f:\Bn\longrightarrow\Cn \textrm { univalent s.t. }  f(0)=0, (df)_0=\id\}. $$

A \textit{(normalized) Loewner chain} is a continuous family of univalent mappings $(f_t:\Bn\longrightarrow\Cn)_{t\geq0}$  with $\Omega_s:=f_s(\Bn)\subseteq f_t(\Bn)$ for $s\leq t$ and such that $f_t(0) = 0$ and $(df)_0=e^t\id$ for all $t\geq 0$.
We also define the \textit{Loewner range} of a normalized Loewner chain as
\[
R(f_t)=\bigcup_{t\geq0}\Omega_t\subseteq\Cn.
\]
The Loewner range is always biholomorphic to $\Cn$,  although, for $n>1$, it can be strictly contain in $\Cn$ (see \cite{ABW}). Furthermore a normalized Loewner chain $(f_t)_{t\geq0}$ is \textit{normal} if $\{e^{-t}f_t(\cdot)\}_{t\geq0}$ is a normal family.

We say that $f\in\S$ \textit{embeds into a Loewner chain} $(f_t)_{t\geq0}$ if $f_0=f$. Then we can define 
$$\S^1:=\{f\in\S: f \ \textrm{embeds into a normalized Loewner chain}\}$$
and, following  ideas of I. Graham, H. Hamada and G. Khor in \cite{GHK} and \cite{S0_2}, we can define the class
$$\S^0:=\{f\in\S: f \ \textrm{embeds into a normal Loewner chain}\}.$$ 
The class $\S^0$ is compact \cite{GK}.
We also denote with 
$$\S(\Cn):=\{\Psi:\Cn\longrightarrow\Cn \textrm { univalent : }  \Psi(0)=0, (d\Psi)_0=\id\}$$
the space of normalized entire univalent functions.
We recall also the following result

\begin{proposition}\cite{GK}\label{p1}
If $(f_t)_{t\geq0}$ is a normalized Loewner chain, then there exist an unique normalized biholomorphism $\Psi:\Cn \longrightarrow R(f_t)$ and normal Loewner chain $(g_t)_{t\geq0}$ such that for each $t\geq0$
$$f_t=\Psi \circ g_t. $$
In particular in the case $t=0$, we have that if $f\in\S^1$ then there exist $\Psi\in\S(\Cn)$ and $g\in\S^0$ such that $f=\Psi\circ g$. In particular, we have the following decomposition 
$$\S^1=\S(\Cn)\circ \S^0. $$
\end{proposition}

\section{Embedding problems in Loewner Theory}

In one dimension all normalized univalent mapping on the disk embed into a normal Loewner chain \cite{PMK} and all normalized Loewner chains are, in fact, normal. Namely,  for $n=1$,
\[\S^0=\S^1=\S. \]

On the other hand, in higher dimension the situation becomes much more complicated. A natural question in Loewner theory in several complex variables, coming from the parallel with dimension one, is the following: 

\medskip
\textit{does every normalized univalent map on the ball embed into a normalized Loewner chain?} 
\medskip

The class $\S^0$ is compact whereas $\S$ and $\S^1$ are not \cite[Chapter 8]{GK}. Hence, there are some normalized univalent functions on the ball that do not embed into a normal Loewner chain (i.e. $\S^0\subsetneq\S$). But, is it possible to embed them into a normalized Loewner chain? The question is still open, and we have the following 

\begin{conjecture}
For $n\geq2$ we have
\[\S=\S^1.\]
By Proposition \ref{p1} this is equivalent to
\[\S=\S(\Cn)\circ\S^0.\]
\end{conjecture}

This conjecture explains very effectively the differences of $\S$ between one and several variables, because,  for $n=1$, 
 $\S(\C)=\{\id_{\C}\}$.
 
\begin{remark}\label{S-S1}
Let $\Psi\in \S(\Cn)$. Then $\Psi|_{\B^n}\in \S^1$. Indeed, $\Psi$ can be embedded into the normalized Loewner chain $\{\Psi(e^t z)\}$.
\end{remark}

A striking difference between one and several complex variables concern  entire univalent mappings: in one dimension, entire univalent maps are only affine transformations while in several variables they are plenty of entire univalent maps which are not affine transformations and not even surjective. This allows for approximation's results which are not possible in dimension $1$. A classical result in this sense is the following

\begin{theorem}\label{thAL}{(And\'ersen-Lempert)}\cite{AL}\\
For $n\geq2$, let $D\subseteq\Cn$ be a starlike domain and let $f:D\longrightarrow \Cn$ be an univalent mapping. Then $f(D)$ is Runge if and only if $f$ can be approximated uniformly on compacta of $D$ by automorphisms of $\Cn$, i.e. for every compact subset $K$ of $D, \epsilon>0$ there exists an automorphism $\Phi$ of $\Cn$ such that $\max_{z\in K}||f(z)-\Phi(z)||<\epsilon$.
\end{theorem}

It is natural to wonder if a similar result is valid for no-Runge mappings, then we have the following

\begin{conjecture}{(Generalized Anderson-Lempert theorem)[GAL]}\label{GAL}\\
For $n\geq2$, let $f\in\S$ with $f(\Bn)$ not Runge. Then there exists a sequence of entire univalent maps which approximate $f$ uniformly on compacta of $\Bn$.
\end{conjecture}

Due to Theorem \ref{thAL}, it is clear that if Conjecture \ref{GAL} holds and $f$ has no Runge image, then entire maps which approximate $f$ uniformly on compacta of $\Bn$ cannot have Runge image. Their images are thus  {\sl Fatou-Bieberbach domains} in $\C^n$ which are not Runge in $\C^n$.

\begin{remark}\label{S1GAL}
In the class $\S^1$, [GAL] holds. Indeed, let $f\in\S^1$, then by Proposition \ref{p1} there exists $\Psi\in\S(\Cn)$ and $g\in\S^0$ such that $f=\Psi\circ g$. Now, $g(\Bn)$ is Runge, then by the And\'ersen-Lempert theorem there exists a sequence of normalized automorphisms $\{\Phi_k \}_{k\in\N}$  that converges to $g$ uniformly on compacta of $\Bn$. Then, also $\{\Psi\circ \Phi_k\}_{k\in\N}\in\S(\Cn)$ converges to $f=\Psi\circ g$ i.e. [GAL] holds in $\S^1$.
\end{remark}

This conjecture is equivalent to the density of $\S^1$ in $\S$. Indeed

\begin{proposition}\label{GALS1}
	For $n\geq2$, the following are equivalent
	
	(1) [GAL] holds;
	
	(2) $\overline{\S^1}=\S$ in the topology of the uniform convergence on compacta.
\end{proposition}
\proof
$(1)\Rightarrow(2)$ Let $f\in\S$. Then there exists a sequence $\{\Psi_k\}\subset\S(\Cn)$ such that $\Psi_k\longrightarrow f$ uniformly on compacta of $\Bn$ (this follows by And\'ersen-Lempert Theorem if $f$ is Runge or by (1) otherwise). But $\Psi_{k|\B^n}\in\S^1$  by Remark \ref{S-S1}, hence (2) holds.

$(2)\Rightarrow(1)$ By Remark \ref{S1GAL} the set
\[\mathcal{A}:=\{f\in\S: f=\Psi_{|\Bn}\ \textrm{for some} \ \Psi\in\S(\Cn) \}  \]
is dense in $\S^1$ i.e. $\S^1\subseteq\overline{\mathcal{A}}$. Therefore by (2)
\[\S=\overline{\S^1}\subseteq\overline{\mathcal{A}}\subseteq\S \ \ \Rightarrow \ \ \overline{\mathcal{A}}=\S. \]
This clearly implies [GAL].
\endproof

We denote

$$\S_R:=\{f\in\S: f(\Bn) \ \textrm{is Runge} \}. $$

In order to find conditions for embedding of univalent functions, a lot of results have been proved during these years. Some of these concern the  regularity of the boundary of $f(\Bn)$. Arosio, Bracci and Wold in 2013 in \cite{ABW} proved  the following 

\begin{theorem}\label{ABW}
Let $n\geq2$ and $f\in\S_R$. If $\Omega:=f(\Bn)$ is a bounded strongly pseudoconvex domain with $\mathcal{C}^\infty$ boundary and $\overline{\Omega}$ is polynomially convex, then $f\in\S^1$.
\end{theorem}

A simple corollary of this theorem is that all functions in $\S_R$ that extend in a univalent and Runge way to a neighborhood of $\overline{\Bn}$ embed into a normalized Loewner chain.

\begin{corollary}
Let $n\geq2$ and $f\in\S_R$. If there exist $r>1$ and $F\in\S_R(r\Bn)$ s.t. $F_{|\Bn}=f$, then $f\in\S^1$.
\end{corollary}
\proof
We use Theorem \ref{ABW}. Obviously $\Omega:=f(\Bn)$ has $\mathcal{C}^\infty$ boundary and $\Omega$ is a strongly pseudoconvex domain (because $\Bn$ is and since the biholomorpishm is defined in a neighborhood of $\overline{\Bn}$ preserves the strongly pseudoconvexity). Now $\forall \ s\in(0,r)$ $\Omega_s:=F(s\Bn)$ are Runge, thus $\overline{\Omega}_s$ has a fundamental system of open neighborhoods which are pseudoconvex and Runge, hence $\overline{\Omega}_s$ is polynomially convex.
\endproof

Following the ideas contained in these results, a natural question is the following: 

\medskip

\textit{does every normalized univalent function of the ball that extends in a univalent way to a neighborhood of the closed ball, embeds into a normalized Loewner chain (without requiring that $f\in\S_R$)? }

\medskip

Therefore we have the following 

\begin{conjecture}{[EXT]}\\
Let $n\geq2$ and $f\in\S$. If exist $r>1$ and $F\in\S(r\Bn)$ s.t. $F_{|\Bn}=f$, then $f\in\S^1$.
\end{conjecture}

\begin{remark}
The conjecture [EXT] is equivalent to requiring that for every $f\in\S$ and $r\in(0,1)$ the mappings $f_r(z):=\frac{1}{r}f(rz)\in\S$ embed into a normalized Loewner chain.  
\end{remark}

Also this conjecture turns out to be equivalent to [GAL].

\begin{proposition}\label{GAL-EXT}
For $n\geq2$, the following are equivalent
	
(1) [GAL] holds (equivalently $\overline{\S^1}=\S$);
	
(2) [EXT] holds.
\end{proposition}

\proof
$(1)\Rightarrow(2)$ Given $f\in\S$ such that there exist $r>1$ and $F\in\S(r\Bn)$ with $F_{|\Bn}=f$, by [GAL] $F$ can be approximated by $\Psi_k\in\S(\Cn)$ uniformly on compacta of $\Bn$. Now there exists $k_0>0$ such that for each $k>k_0$ we have 
\[f(\overline{\Bn})=:\overline{\Omega}\subset\Psi_k(\Bn)\] 
and $\Psi^{-1}_k\longrightarrow F^{-1}$ uniformly on $\overline{\Omega}$. But $F^{-1}(\overline{\Omega})=\overline{\Bn}$ is a strongly convex set and since strong convexity is an open condition (in the $\mathcal{C}^2$ topology), there exist $k_1>k_0$ and  $\Psi:=\Psi_{k_1}$ such that $\Psi^{-1}(\overline{\Omega})$ is strongly convex. Therefore $f$ can be embedded into the Loewner chain 
\[f_t(z):=\Psi(e^t\cdot\Psi^{-1}(f(z))).\]

$(2)\Rightarrow(1)$ Given $\{r_k\}_{k\in\N}\subset(0,1)$ a sequence that converges to $1$, for each $f\in\S$ we can define $f_k(z):=\frac{1}{r_k}f(r_kz)\in\S$. Obviously $f_k \longrightarrow f$ uniformly on compacta of $\Bn$, and $f_k\in\S^1$ for each $k\in\N$ by [EXT]. Hence $\overline{\S^1}=\S$.
\endproof

We recall the following result, that easily descends from the Docquier-Grauert Theorem \cite{DG}.

\begin{proposition}\label{thDG}
Let $(f_t)_{t\geq0}$ be a normalized Loewner chain. Then for each $0<s\leq t$ the couple $(f_s(\Bn),f_t(\Bn))$ is Runge. Therefore  for each $ t\geq0$ also the couple $(f_t(\Bn),R(f_t))$ is Runge.
\end{proposition}

Given $(f_t)_{t\geq0}$ a normalized Loewner chain, by the previous proposition we have that $f_0(\Bn)$ is Runge into its Loewner range $R(f_t)$ that is biholomrphic to $\Cn$. Thus a necessary condition in order to obtain embedding of all the univalent functions into a Loewner chain is that all the open sets of $\Cn$ biholomorphic to $\Bn$ have to be Runge in $\Cn$ or in some Fatou-Bierberbach domain. Thus we have the following 

\begin{conjecture}{[FBR]}\\
Let $D\subseteq\Cn$ a domain biholomorphic to the ball $\Bn$, which is not Runge in $\Cn$. Then there exists $\Omega\subseteq\Cn$ a Fatou-Bieberbach domain with $D\subseteq\Omega$ such that $(D,\Omega)$ is a Runge pair.
\end{conjecture}

Related to the previous conjecture, we have the following strong version of [GAL].

\begin{conjecture}{(Strong Generalized And\'ersen-Lempert Theorem)[GAL$^s$]}\\
For $n\geq2$, let $f\in\S$ with $f(\Bn)$ not Runge. Then there exist a Fatou-Bieberbach domain $\Omega$ and a sequence of univalent mappings $\Psi_k\in\S(\Cn)$ with $\Psi_k(\Cn)=\Omega$ for each $k$, which converges to $f$ uniformly on compacta of $\Bn$.
\end{conjecture}

We have also the following weaker formulations

\begin{conjecture}{[FBR$_a$]}\\
Let $f\in\S$, then for each $r\in(0,1)$ exists $\Omega_r\subseteq\Cn$ a domain biholomorphic to $\Cn$ with $f(r\Bn)\subseteq\Omega_r$ such that $(f(r\Bn),\Omega_r)$ is a Runge pair.
\end{conjecture}

\begin{conjecture}{[GAL$_a^s$]}\\
Let $f\in\S$ with $f(\Bn)$ not Runge, then for each $r\in(0,1)$ there exist a Fatou-Bieberbach domain $\Omega_r$ and a sequence of univalent mappings $\Psi_k^{(r)}\in\S(\Cn)$ with $\Psi_k^{(r)}(\Cn)=\Omega_r$ for each $k$, which converges to $f_r$ uniformly on compacta of $\Bn$.
\end{conjecture}

Obviously [FBR] and [GAL$^s$] imply respectively [FBR$_a$] and [GAL$_a^s$], and furthermore they are equivalent two by two.

\begin{proposition}\label{FBR-GALs}
For $n\geq2$, we have
	
(1)[FBR] and [GAL$^s$] are equivalent.
	
(2)[FBR$_a$] and [GAL$_a^s$] are equivalent.
\end{proposition}
\proof
The two proofs are essentially the same, so we only prove (1).

([FBR]$\Rightarrow$[GAL$^s$]) 
Let $f\in\S$ and $\Psi\in\S(\Cn)$ be a Fatou Bierberbach mapping such that ($f(\Bn),\Psi(\Cn)$) is Runge. Now since Runge-ness is a property invariant under biholomorphism, ($(\Psi^{-1}\circ f)(\Bn),\Cn$) is Runge. By the And\'ersen-Lempert theorem $\Psi^{-1}\circ f$ is approximable by automorphisms $\Phi_k$, then $\Psi\circ \Phi_k\in\S(\Cn)$ converges to $f$, i.e. [GAL$^s$] holds.

([GAL$^s$]$\Rightarrow$[FBR]) 
Let $f\in\S$ and $\{\Psi_k\}_{k\in\N}\subset\S(\Cn)$ a sequence of normalized entire univalent mappings with the same image $\Omega$ that approximate $f$. Now $f(\Bn)\subseteq\Omega$, then the sequence $\{\Psi^{-1}_0\circ\Psi_k\}_{k\in\N}$ converges to $\Psi_0^{-1}\circ f$. Note that $\Psi^{-1}_0\circ\Psi_k$ are normalized automorphisms for each $k\in\N$, then by the necessity condition of Andérsen-Lempert theorem $\Psi_0^{-1}\circ f$ has to be Runge. Finally, by invariance of Runge-ness under biholomorphisms, $(f(\Bn),\Psi_0(\Cn)=\Omega)$ is Runge.
\endproof

\begin{proposition}\label{pf}
	For $n\geq2$, we have
	
	(1)[GAL$_a^s$]  implies [GAL].
	
	(2)[EXT] implies [FBR$_a$].
	
	Moreover, by Proposition \ref{GAL-EXT} and Proposition \ref{FBR-GALs}  [FBR$_a$], [GAL], [GAL$_a^s$] and [EXT] conjectures are equivalent to $\overline{\S^1}=\S$.
\end{proposition}
\proof
([GAL$_a^s$]$\Rightarrow$[GAL]) 
Let $f\in\S$ with $f(\Bn)$ non Runge and $K\subseteq\Bn$ a compact. Then there exists $r\in(0,1)$ such that $K\subseteq r\Bn$. Thanks to [GAL$_a^s$], there exists a normalized Fatou Bieberbach mapping $\Psi\in\S(\Cn)$ such that
\[\max_{z\in \frac{1}{r}K}|| f_r(z)-\Psi(z)||<\epsilon/r\]
then
\[\max_{z\in K}||f(z)-r\Psi(z/r)||=r\max_{z\in \frac{1}{r}K}||f(rz)/r-\Psi(z)||<\epsilon\]
i.e. [GAL] holds.

([EXT]$\Rightarrow$[FBR$_a$])
Fix $f\in\S$. For each $r\in(0,1)$ by [EXT] $f_r$ embeds into a normalized Loewner chain, then $f_r(\Bn)$ is Runge in its Loewner range.
\endproof

To summarize, we have the following scheme

\[
\begin{tikzcd}[column sep=3pc,row sep=3pc]
\textrm{[GAL]} \arrow[r,Leftrightarrow,"\ref{GAL-EXT}"] \arrow[d,Leftarrow,"\ref{pf}"] & \textrm{[EXT]} \arrow[d,Rightarrow,"\ref{pf}"] \\
\textrm{[GAL}^s_a\textrm{]} \arrow[r,Leftrightarrow,"\ref{FBR-GALs}"] \arrow[d,Leftarrow]  & \textrm{[FBR}_a\textrm{]} \arrow[d,Leftarrow]\\
\textrm{[GAL}^s\textrm{]} \arrow[r,Leftrightarrow,"\ref{FBR-GALs}"] & \textrm{[FBR]}
\end{tikzcd}
\]

\begin{remark}
(1) $\S_R$ is close and by the And\'ersen-Lempert theorem $\S_R\subseteq\overline{\S^1}$. Therefore, if $\S^1$ were closed then $\S_R\subseteq\S^1$ i.e. every univalent function with Runge image would be embedded into a normalized Loewner chain.\\
(2) If [FBR] holds and $\S_R\subseteq\S^1$, then every normalized univalent function would be embedded into a normalized Loewner chain. Indeed if $f\in\S$ with $f(\Bn)$ not Runge, by [FBR] there exists $\Psi\in\S(\Cn)$ such that $(\Psi^{-1}\circ f)(\Bn)$ is Runge, hence, if $\S_R \subseteq \S^1 \ \Psi^{-1}\circ f\in\S^1$. Then $f\in\S^1$.
\end{remark}

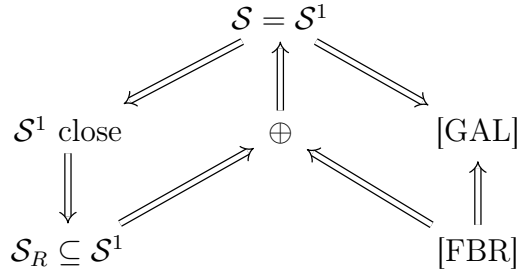
\[
\begin{tikzcd}[column sep=2.5pc,row sep=2pc]
{}      & \S=\S^1 \arrow[ld,Rightarrow] \arrow[rd,Rightarrow]& \\
\S^1 \ \textrm{close} \arrow[d,Rightarrow]   &  \oplus \arrow[u,Rightarrow]\arrow[ld,Leftarrow]\arrow[rd,Leftarrow] & \textrm{[GAL]} \arrow[d,Leftarrow]\\
\S_R\subseteq \S^1    &   & \textrm{[FBR]}
\end{tikzcd}
\]

\section{Entire-convexshapelike domains}
In this section we give a proof of the following theorem, that is the generalization of Theorem \ref{ABW} under the additioned hypothesis that [GAL] holds.

\begin{theorem}\label{SPD}
Let $n\geq2$ and suppose [GAL] holds. If $f\in\S$ has the property that $f(\Bn)$ is a bounded strongly pseudoconvex domain with $\mathcal{C}^\infty$ boundary then $f\in\S^1$.
\end{theorem}

In \cite{ABW}, in order to prove Theorem \ref{ABW} it was introduced the concept of \textit{convexshapelike domains}: $\Omega\subseteq\Cn$ is convexshapelike if there exists $\Phi$ an automorphism of $\Cn$ such that $\Phi(\Omega)$ is convex in $\Cn$. This kind of domains are very useful in the study of embedding problems, indeed given $f\in\S$ if $f(\Bn)$ is convexshapelike then $f$ embeds into the Loewner chain
\[f_t(z):=\Phi^{-1}(e^t\cdot\Phi(f(z))).  \]
According to \cite{ABW}, a Loewner chain of this form is called a {\sl filtering Loewner chain}.
\begin{definition}
A normalized Loewner chain $(f_t)_{t\geq0}$  in $\Bn$  is a \textit{filtering normalized Loewner chain} if $\Omega_t:=f_t(\Bn)$ has the following properties:
	
$(1) \ \bar{\Omega}_s\subseteq\Omega_t$ for each $0\leq s<t$;
	
$(2)$ for each open set $U$ containing $\Omega_s$ there exists  $ t_0>s$ such that $\Omega_t\subset U$ for all $t\in(s,t_0)$. 
\end{definition}

A natural generalization of the concept of convexshapelike domains is the following

\begin{definition}
$\Omega\subseteq\Cn$ is \textit{entire-convexshapelike domains} if there exists $\Psi:\Cn\longrightarrow\Cn$ an entire univalent mapping with $\Omega\subseteq\Psi(\Cn)$ s.t. $\Psi^{-1}(\Omega)$ is convex.
\end{definition}

As before, if $f\in\S$ and $f(\Bn)$ is entire-convexshapelike then $f$ embeds into the filtering normalized Loewner chain 
\[f_t(z):=\Psi(e^t\cdot\Psi^{-1}(f(z))),\]
where $\Psi:\Cn\longrightarrow\Cn$ is a normalized entire univalent mapping such that $\Psi^{-1}(f(\Bn))$ is convex.\\

\begin{lemma}\label{lemapr}\cite{ABW}\\
Let $\Omega\subseteq\Cn$ be a bounded strongly  pseudoconvex domain with $\mathcal{C}^\infty$ boundary which is biholomorphic to $\Bn$. Then any $f\in \mathcal{C}^2(\overline{\Omega})\cap$Hol$(\Omega,\Cn)$ can be approximated uniformly on $\overline{\Omega}$ in $\mathcal{C}^2$-norm, by functions in Hol$(\overline{\Omega},\Cn)$. 
\end{lemma}

Therefore, in order to obtain Theorem \ref{SPD}, we prove the following

\begin{proposition}\label{proec}
Let $n\geq2$ and suppose [GAL] holds. If $\Omega$ is a bounded strongly pseudoconvex domain with $\mathcal{C}^\infty$ boundary, then $\Omega$ is entire-convexshapelike.
\end{proposition}

\proof 
By Fefferman's Theorem \cite{FEF} $f$ extends to a diffeomorphism $f:\overline{\Bn}\longrightarrow \overline{\Omega}$. By Lemma \ref{lemapr} $f^{-1}$ can be approximated in $\mathcal{C}^2$ norms uniformly on $\overline{\Omega}$ by holomorphic  maps defined on neighborhoods of $\overline{\Omega}$. Then there exists an open neighborhood $U$ of $\overline{\Omega}$ and $h:U\longrightarrow\Cn$ an univalent mapping s.t. $h(\Omega)$ is a smooth strongly convex domain (because strongly convexity is an open condition). Without loss of generality, we can choose $h(U)=r\Bn$ with $r>1$. Now, by [GAL], $h^{-1}$ can be approximated uniformly on compacta of $h(U)=r\Bn$ by $\Psi_k\in\S(\Cn)$, then there exists $k_0>0$ such that for each $k>k_0$ we have 
\[h^{-1}(\overline{\Bn})=\overline{\Omega}\subset\Psi_k(\Bn)\] 
and $\Psi^{-1}_k\longrightarrow h$ uniformly on $\overline{\Omega}$. But $h(\overline{\Omega})=\overline{\Bn}$ is a strongly convex set, therefore there exist $k_1>k_0$ and  $\Psi:=\Psi_{k_1}$ such that $\Psi^{-1}(\overline{\Omega})$ is strongly convex, i.e. $\Omega$ is entire-convexshapelike.
\endproof

We recall a classical result of several complex variables.

\begin{lemma}{(Narasimhan)}\\
Let $\Omega\subset\Cn$ be a domain with boundary $\mathcal{C}^2$ in $p\in\partial\Omega$, and suppose that $\Omega$ is strongly Levi pseudoconvex in $p$. Then there exists an open neighborhood $U$ of $p$ and a biholomorphism $f:U\longrightarrow V\subset\Cn$ such that $f(U \cap\Omega)$ is a convex domain.
\end{lemma}

Proposition \ref{proec} can be seen as a global version of Narasimhan's lemma (for domains biholomorphic to a ball and with $\mathcal{C}^\infty$ boundary).

We conclude with the following

\begin{remark}
Suppose the statement of Proposition \ref{proec} holds (i.e. every domain biholomorphic to a ball, with $\mathcal{C}^{\infty}$ boundary and strongly pseudoconvex, is entire-convexshapelike), then $\S^1$ is dense in $\S$. Indeed [EXT] holds: fix $f\in\S$, then for each $r\in(0,1)$ $f_r(z):=\frac{1}{r}f(rz)$ is entire-convexshapelike and therefore it is in $\S^1$.

\end{remark}

\textbf{Acknowledgements:}

The author would like to express his gratitude to Prof. Filippo Bracci for his availability and for introducing him into Loewner theory.

\end{document}